\documentclass[authoryear, preprint]{elsarticle}
\usepackage[margin=3cm]{geometry}
\usepackage{amsmath, amsthm, amssymb}
\usepackage{graphicx}

\newcommand{\tr}{^{\prime}}

\def\b#1{\mbox{\boldmath $#1$}}    
\def\bl#1{\mbox{\footnotesize \boldmath {$#1$}}} 
\def\m#1{\mbox{#1}}                
\def\cg#1{\ensuremath{\mathcal{#1}}}      

\newcommand{\diag}{{\rm diag}}    

\newtheorem{definition}{Definition}

\renewcommand{\mid}{\,|\,}

\usepackage{color}
\definecolor{dg}{rgb}{0.1,0.5,0.1}
\title{Testing order restrictions in contingency tables}

\author[rc]{R.~Colombi}
\author[af]{A.~Forcina}
\address[rc]{Dipartimento di Ingegneria, University of Bergamo, Italy}
\address[af]{Dipartimento di Economia, Finanza e Statistica,
University of Perugia, Italy}

\begin{document}
\begin{abstract}
Several interesting models for contingency tables are defined by a
system of equality and inequality constraints on a suitable set of
marginal log-linear parameters. After reviewing the most common
difficulties which are intrinsic to order restricted testing
problems, we propose two new families of testing procedures, based
on similar attempts appeared in the econometric literature, in order
to increase the probability of detecting several relevant violations
of the supposed order relations. One set of procedures is based on
the decomposition of the log-likelihood ratio when testing the given
set of inequalities and the nested model derived by forcing
inequalities into strict equalities. The other set uses the
asymptotic joint normal distribution of the estimates of the
marginal log-linear parameters to be constrained.
\end{abstract}

\begin{keyword}
Stochastic orderings \sep chi-bar squared distribution \sep positive
association
\end{keyword}
\maketitle
\section{Introduction}
Models for contingency tables involving order restrictions may arise
in several contexts. For two-way tables, when both variables have
ordered categories, \citet{DarFor98} extended the approach of
\citet{DyKoRo} for testing the set of inequalities implied by the
assumption that the conditional distributions by row satisfy
suitable stochastic orderings and  \citet{Coca} considered testing
stochastic ordering by row and columns simultaneously. More
generally, inequality constraints arise when we assume that a pair
of ordered categorical variables are positively associated
\citep[][]{BaCoFo} or that the strength of the association increases
with respect to a third variable \citep{ColFor:01}. Order
restrictions may also be implied by graphical models containing
latent variables when we restrict attention to the marginal
distribution of the observed variables, for a general treatment see
\citet{Evans}. In the very special case of item response models,
\citet{BartForc} considered testing the order relation known as MTP2
which is implied by the assumptions of conditional independence and
monotonicity. For a very extensive review of the literature on the
subject see \citet{AgriCol}.

The problem of testing a model defined by a system of linear inequality constraints imposed on a set of marginal log-linear parameters against the saturated model is a difficult one because, in the most common formulations, the null hypothesis states that the true parameter value lies into a convex cone and {may be on the boundary} \citep[][]{Shapiro85}. The problem of determining a least favorable distribution within the null may be a difficult one \citep[See][Sec. 4.8.5]{SiSe} because this does not necessarily coincides with the model where all inequalities hold as strict equalities. However, the more restrictive hypothesis where all inequalities hold as strict equalities may be of interest on its own, for instance, the assumption of positive quadrant dependence in a two-way table \citep{BaDaFo}, defined by the constraint that all the log-odds ratios of type global are non negative, corresponds to a convex cone whose vertex is the model of independence. The problem of testing that a system of linear inequalities may be replaced by equalities against the assumption that at least one inequality is strict \citep[see for example][Sec. 4.3]{SiSe} is closely related to the problem of testing the corresponding set of inequalities against the saturated model in the sense that  the likelihood ratio statistics for testing a system of equalities against the unrestricted model may be partitioned into two components, one for testing equalities against inequalities and the other for testing inequalities against the saturated model, whose asymptotic distributions are  mixture of chi squared random variables.

{It is well known that testing a set of equalities against a
corresponding set of inequalities, may lead to  reject the null even
if the inequality constraints are clearly violated in the data }
\citep[see][3.3]{AgriCol}. This fact has been used by \citet{KatAgr}
as an argument against using frequentist approaches in this context.
The problem has been studied in detail within the econometric
literature and valid frequentist solutions have been proposed. For
instance, in the context of testing Lorenz curve orderings,
\citet{DaFoLor} designed a procedure based on the joint distribution
of the two components of the likelihood ratio statistics mentioned
above. {Starting from similar concerns}, \citet{Bennet} proposed a
procedure based on a multiple comparisons approach. In this paper we
propose an extension of both procedures and compare their merits
with a simulation. It turns out that the procedures based on the
joint distribution of the likelihood ratio statistics performs much
better against most alternatives.

In section two we introduce the notation and  present a formal
statement of the general problem. In section 3 we review the basic
features of the problem of testing a set of inequality constraints
and of  {the} chi-bar squared distribution. In section four we present two new sets of testing procedures, one based on  the two likelihood
ratio statistics for testing equalities against inequalities and
inequalities against the saturated model and the other on Multiple
Comparisons approaches. Two real data sets are analyzed in Section 4  to illustrate the
procedures. The simulation study of section five provides clear
evidence that the proposed approaches can provide satisfactory
solutions, though the procedures based on the Likelihood ratio
statistics perform substantially better.

\section{Notation and preliminary results}
Consider a contingency table determined by the joint distribution of
$d$ discrete random variables and let $\b i$ be a vector of $d$
indices; $p_{\bl i}$ denote the probability that an observation
falls in cell $\b i$ and $\b p$ the vector whose elements $p_{\bl
i}$ are arranged in lexicographic order relative to $\b i$. We
assume that the elements of $\b p$ are strictly positive. Let
$\b\eta$ be a vector of marginal log-linear parameters as defined in
\citet{BaCoFo} and we assume that the mapping between $\b\eta$ and
$\b p$ is a diffomorphism. Any such vector may be defined as
$\b\eta$ = $\b C\log(\b M\b p)$ where $\b M$ is a matrix of 0's and
1's which produce the appropriate marginal or aggregated
probabilities and $\b C$ is a matrix of row contrasts. This
formulation allow to consider ordinary log-linear parameters as well
as logits and higher order interactions  of type global or
continuation. Assume that $\b\eta$ is contained in an open subset of
$\Re^{t-1}$, where $t$ is the number of cells of the table.

{Let $\cg C$ denote a closed convex cone and assume that $\cg
L_0,\:\cg L_1$ are, respectively, the linear space of largest and
smallest dimensions such that $\cg L_0\subset \cg C \subset \cg
L_1$. Let dim$(\cg L_0)=q$, dim$(\cg L_1)=r$; finally let $\cg S
\supseteq \cg L_1$ be the parameter space of the saturated model.}
Define
$$
H_0: \b\eta\in\cg L_0,\quad  H_1:\b\eta \in \cg C,\quad H_2:
\b\eta\in \cg S
$$
and note that $H_1$ will usually be defined by a set of equality and
inequality constraints while $H_0$ is obtained by turning all
inequalities into strict equalities.
\section{Hypotheses testing}
Let $L_{01}$ be the log-likelihood ratio for testing $H_0$ against $H_1$
and $L_{12}$ the log-likelihood ratio for testing $H_1$ against $H_2$.
It is well known that the asymptotic distribution of
$L_{02}=L_{01}+L_{12}$ is $\chi^2_s$, where $s=t-q-1$. Let $\b F_0$
be the expected information matrix of $\b\eta$ under $H_0$ and $\b
V_0$ =$\b F_0^{-1}$. It can also be shown \citep[see for
example][4.3]{SiSe} that
\begin{equation}
Pr(L_{01}>c\mid H_0) = \sum_q^{r} w_j(\b V_0,\cg C)
Pr(\chi_{j-q}^2>c), \label{eq:chibar}
\end{equation}
where $\chi^2_j$ denotes a chi square random variable with $j$
degrees of freedom. The above distribution, known as chi-bar
squared, depends on the probability weights $w_j(\b V_0,\cg C)$
whose definition and computation we discuss below. It can be shown
that, under $H_0$, $L_{12}$ is asymptotically  distributed like the
sum of a $\chi^2_{t-1-r}$ and a chi-bar squared with weights $w_j(\b
V_0,\cg C)$ used in the reverse order \citep[see for
instance][4.8.5,]{SiSe}. Let $\hat{\b\eta}$ denote the unrestricted
maximum likelihood estimate; an interesting geometric interpretation
is that $L_{01}$ and $L_{12}$ are asymptotically equivalent to the
squared norm of the projection of $\hat{\b\eta}$ onto, respectively,
the convex cone defined by $H_1$ and onto its dual, in the metric
defined by $\b F_0^{-1}$. From this, the following  simple
expression for the asymptotic joint distribution of $L_{01}$ and
$L_{12} $ can a be derived \citep[see for example][Lemma
2]{DaFoLor}:
\begin{equation}
Pr(L_{01} \leq c_1,L_{12}\leq c_2\mid H_0) = \sum_q^{r} w_j(\b
V_0,\cg C) Pr(\chi_{j-q}^2 \leq c_1)Pr(\chi_{t-j-1}^2 \leq c_2)
\label{joint}
\end{equation}
The previous joint distribution will play a key role in the next
section where we consider testing procedures which use the
statistics $L_{01}$ and $L_{12}$ simultaneously.
\subsection{The probability weights}
Fast and accurate computation of the weights $w_j(\b V_0,\cg C)$ is
crucial to the use of the testing procedures of section 4.1
below, so we discuss two alternative methods that can be used to
evaluate these probabilities. For simplicity, we restrict to
the context where $r=t-1$, meaning that $H_1$ is defined only by
inequality constraints. The more general case where equality
constraints are also present, so that $\cg L_1\subset \cg S$, can be
reduced by replacing $\b y$ with its projection onto $\cg L_1$ and
$t-1$ with $r$, see \citet{Shapiro88} for a detailed treatment.

If we are prepared to assume that the distribution of $L_{01}$ and
$L_{12}$ should be determined under $H_0$, the probability weights
may be defined as follows. Let $\b x$ be distributed as multivariate
normal $\cg N(\b 0, \b V_0)$, then $ w_j(\b V_0,\cg C)$ is the
probability that the projection of $\b x$ onto $\cg C$ falls on a
face which spans a linear space of dimension $j$ \citep[see for
example][Prop. 3.6.1]{SiSe}. Recent results \citep{GeBr} on
multivariate normal integrals, make computation of exact weights
fast and accurate  for moderate values of $r$. An algorithm for
computing the probability weights derived from \citep[][pag
409-416]{Kudo} is implemented in the R-packages {\bf ic-infer}
\citep{Ulrike} and {\bf hmmm} \citep{rsm}; a slightly different
implementation is outlined in the Appendix. However when the number
of inequalities is, say, larger than 15, exact estimation becomes
very hard. For this reason, \citet{DarFor98} suggested a simulation
procedure by which a reasonably large number of sample points from
the appropriate normal distribution are projected onto $\cg C$; then
$ w_j(\b V_0,\cg C)$ is estimated by the proportion of sample points
falling on any face of $\cg C$ of dimension $j$.

The following method, which extends a similar one outlined by
\citet{DarFor98}, can be used  to determine the minimum number of
sample points required for accurate estimation of probability
weights. Let $\b w$ be the vector with elements $ w_j(\b V_0,\cg C)$
and $\hat{\b w}$ its estimate. Let  $\b c$ be the vector with
elements $P(\chi_j^2\leq c_1) P(\chi_{q-j}^2\leq c_2)$, where
$c_1.\:c_2$ are determined as in the tunable LR testing procedure of
Definition 3 below, and let $v$ = $\b c\tr[\diag(\b w)-\b w\b
w\tr]\b c$. Suppose we require that the estimation error $\b
c\tr(\hat{\b w}-\b w)$ of estimating the probability of accepting
$H_0$ when true must satisfy the condition
$$
P(\mid \b c\tr(\hat{\b w}-\b w)\mid \leq \epsilon)\leq 1-2\delta.
$$
By using the multivariate normal approximation to the
multinomial and elementary calculations, we have that the number of
points to be projected  should not be smaller than
$v(z_{1-2\delta}/\epsilon)^2$ where $z_{1-2\delta}$ is the
$1-2\delta$ percentile of the standard normal distribution.

In general \citep[see][4.3.1]{SiSe}  nuisance parameters are going
to affect  $\b V_0$ and thus the probability weights. The formally
correct procedure would be to search for the least favorable null
distribution, a task which, however, may be very hard. In addition,
it often turns out that the values of the nuisance parameters that
produce the least favorable distribution are very extreme and
substantially different from any plausible estimate obtained from
the data. An alternative procedure would be to compute the null
distribution after replacing the nuisance parameters with their
maximum likelihood estimate; a simulation study
\citep[][4.5]{DarFor98} indicates that the $p$-value computed in
this way is sufficiently close to the one computed at the true value
of the nuisance parameters when the sample size is moderately large.

In certain contexts one could remove the dependence on nuisance
parameters by conditioning; for instance, if we are interested in
the dependence structure of a two way table, we might condition to
the row and columns totals, an approach explored by \citet{BaDaFo}. The resulting null distributions are, however, hard to handle even with the power of modern computers and can be applied for small sample sizes and in specific contexts.

\section{Testing procedures}
{In this section we present two different families of testing
procedures which may lead to one of the following decisions: (i)
accept $H_0$, (ii) reject $H_0$ in the direction of $H_1$, meaning
that there is strong evidence to support the assumed inequality
constraints, (iii) reject $H_0$ in the direction of $H_2$, if there
is convincing evidence that the inequality constraints are violated.
As mentioned in the introduction, if we considered the problem of
testing $H_1$ against the saturated model on its own, then we would
face the additional problem of having to determine the least
favorable distribution within $H_1$: \citet{Wolak} has shown that,
with non linear problems, the least favorable distribution may not
coincide with $H_0$. In this paper we avoid this complication
because  the two testing problems ($H_0$ against $H_1$ and $H_1$
against $H_2$) are combined into a single testing problem, as
detailed in the sub-sections below.}
\subsection{Likelihood ratios}
As emphasized, for instance, by \citet{AgriCol}, the evidence in
favour of $H_1$ provided by $L_{01}$  may be highly misleading
because a large value of this statistic, which would lead to reject
$H_0$ in favour of $H_1$, is compatible with substantial violations
of $H_1$. A geometric explanation is that this event will happen
whenever $\hat{\b\eta}$ is far away from $\cg L_0$ and is not
contained neither in $\cg C$, nor into its dual. Thus, before
deciding that the assumed set of inequality constraints are
satisfied, one should also examine $L_{12}$: a large value of this
statistic provides evidence that $H_1$ is violated in the direction
of $H_2$.

Consider first the following procedure:
\begin{definition} The naive procedure:
\begin{enumerate}
\item accept $H_0$ if $L_{01}\leq c_1$  where $Pr(L_{01}\leq c_1\mid
H_0)$ = $1-\alpha_1-\alpha_2$,
\item reject $H_0$ in favour of $H_1$ if $L_{01}>c_1$ and $L_{12}\leq
c_{2}$, where $Pr(L_{01}>c_1, \:L_{2}\leq c_{2}\mid H_0)$ =
$\alpha_1$,
\item otherwise reject $H_0$ for $H_2$.
\end{enumerate}
\end{definition}
This procedure is closely related to a widely adopted standard
approach  that rejects $H_0$ towards $H_1$ for large values of
$L_{01}$ and,  when $H_0$ is rejected, uses the statistics $L_{12}$
to test $H_1$ against $H_2$. However, as far we know, the
implementations of this approach, considered in the literature, do
not use the joint distribution (\ref{joint}) of the two statistics to control the
error rate of a false rejection of $H_1$.

An alternative approach, also related to a standard use of $L_{01}$
and $L_{12}$, is:
\begin{definition} The basic procedure: let $c_2$ be such that
$Pr(L_{12}> c_2\mid H_0)$ = $\alpha_2$ then
\begin{enumerate}
\item accept $H_0$ if $L_{12}\leq c_2$ and $L_{01}\leq c_1$,
where $Pr(L_{01}\leq c_1, L_{12}\leq c_2\mid H_0)$ = $1-\alpha_1-\alpha_2$,
\item reject $H_0$ towards $H_1$ if $L_{01}>c_1$ and  $L_{12}\leq c_2$,
\item reject $H_0$ towards $H_2$ if $L_{12}>c_2$.
\end{enumerate}
\end{definition}

Though the above procedure allows a direct control of the error
rates towards $H_1$ and $H_2$ and does not suffer from the
limitations described by \citet{AgriCol}, its power in detecting
$H_2$, when it is true, may still be too low. \citet{DaFoLor}, in
the context of comparing economic inequality, suggested a procedure
where the amount of protection against rejecting $H_0$ in favour of
$H_1$ when $H_2$ holds may be tuned to the specific context. The
procedure we describe below is a revised version of their procedure.
\begin{definition} The tunable procedure: let
$0 \leq \alpha_{12} \leq \alpha_2$ and let $c_2$ be such that
$Pr(L_{12}> c_2\mid H_0)$ = $\alpha_2-\alpha_{12}$, then
\begin{enumerate}
\item accept $H_0$ if $L_{12}\leq c_2$ and $L_{01}\leq c_1$,
where $Pr(L_{01}\leq c_1, L_{12}\leq c_2\mid H_0)$ = $1-\alpha_1-\alpha_2$,
\item reject $H_0$ in favour of $H_1$ if $L_{01}>c_1$ and
$L_{12}\leq c_{12}$, where $Pr(L_{01}>c_1, \:L_{12}\leq c_{12}\mid
H_0)$ = $\alpha_1$,
\item otherwise reject $H_0$ for $H_2$.
\end{enumerate}
\end{definition}
Note that the tunable procedure reduces to the naive procedure when
$\alpha_{12}=\alpha_2$ and is equivalent to the basic procedure when
$\alpha_{12}=0$. The tunable procedure allows to fix the error rates
towards each alternative and, in addition, by increasing the tuning
probability $\alpha_{12}$ it can decrease the probability of
rejecting towards $H_1$ when some of the assumed constraints are
violated. However, larger values of $\alpha_{12}$ correspond also to
larger values of $c_2$, thus, the side effect is an increasing
tendency to accept $H_0$ even when it should be rejected towards
$H_2$, which is the main drawback of the naive procedure. This fact
will be clarified by the simulation study in the last section.
\subsection{Multiple Comparison procedures}
There is a collection of procedures which, due to their
computational simplicity, have received much attention in
Econometric applications about inequality constrained inference
problems, see for instance \citet{Bishop:91}. To keep notation
simple, we restrict attention to the case where $H_1$ is defined by
$\b D\b\eta\geq \b 0$, where the $k\times (t-1)$ matrix $\b D$ is of
full row rank. Let $\b\tau$ = $\b D\b\eta$ and $\hat{\b\tau}$ denote
the unrestricted maximum likelihood estimator. Let $n$ denote the
sample size; under $H_0$ and the usual regularity conditions, it
follows that the asymptotic distribution of $\b z=\surd n
\hat{\b\tau}$
is multivariate normal with covariance matrix  $\b \Sigma_0$ = $\b
D\b V_0\b D\tr$. Let $max(\b z)$ and $min(\b z)$ denote the greatest
and the lowest component respectively of $\b z$. The following
approaches have in common the fact that the different rejection
regions are based on the $max(\b z)$ and $min(\b z)$ statistics.

\begin{definition}
The naive Multiple Comparisons procedure: let $c$ be such that
$Pr(min(\b z)\geq-c, \: max(\b z)\leq c\mid H_0)$ = $1-\alpha$, then
\begin{enumerate}
\item accept $H_0$ if $max(\b z)\leq c$ and $min(\b z) \geq -c$
\item reject $H_0$ for $H_1$ if $max( \b z)>c$ and $min( \b z)\geq -c$,
\item reject $H_0$ for $H_2$ otherwise.
\end{enumerate}
\end{definition}
As noted by \citet{DaFoLor}, this procedure has the  limitation that
the probabilities of rejecting $H_0$ towards $H_1$ and $H_2$ cannot
be controlled. \citet{Bennet}, in the context of comparing income
inequality, has recently proposed a procedure which tries to
overcome this limitation. Though, formally, Bennet's procedure is
based on the one-sided Kolmogorov-Smirnov statistic, it can be
adapted to the present context and formulated in terms of the
asymptotic normal distribution of the vector $\b z$, defined above.
\begin{definition}
Bennet's procedure: let $c_1$ be such that $Pr(min(\b z)\geq-c_1, \:
max(\b z)\leq c_1\mid H_0)$ = $1-\alpha_1-\alpha_2$, then
\begin{enumerate}
\item accept $H_0$ if $max(\b z)\leq c_1$ and $min(\b z) \geq -c_1$,
\item reject $H_0$ towards $H_1$ if max$(\b z)> c_1$ and min$(\b z)>-c_2$
where $Pr(\m{max}(\b z)>c_1,\m{min}(\b z)>-c_2\mid H_0)=\alpha_1$,
\item reject $H_0$ towards $H_2$ otherwise.
\end{enumerate}
\end{definition}
In this procedure $\alpha=\alpha_1+\alpha_2$ is the total
probability of rejecting $H_0$ when true; Bennet sets $\alpha_1$ =
$\alpha\beta$ and $\alpha_2$ = $\alpha(1-\beta)$, with $\beta$ =
$\frac{\alpha_1}{\alpha}$.

We now present a more flexible procedure which, for any preassigned value of $\alpha_1,\:\alpha_2$, like in Definition
3 above, allows to increase the level of protection against stating
that $H_1$ holds when one or more inequalities are
violated in the population.
\begin{definition} { The tunable Multiple Comparisons (MC) procedure:
let  $0\leq \alpha_{12}\leq \alpha_2$ and $c_2$ be such that
$Pr(\m{min}(\b z)< -c_2\mid H_0)$ = $\alpha_2-\alpha_{12}$ and then}
\begin{enumerate}
\item accept $H_0$ if $\m{min}(\b z)\geq -c_2$ and $\m{max}(\b z)
\leq c_1$, where $Pr(\m{max}(\b z)\leq c_1,\: \m{min}(\b z)
\geq -c_2\mid H_0)$ = $1-\alpha_1-\alpha_2$,
\item reject $H_0$ in favour of $H_1$ if $\m{max}(\b z)>c_1$ and
$\m{min}(\b z)\geq -c_{12}$, where $Pr(\m{max}(\b z)>c_1, \:\m{min}
(\b z)\geq -c_{12}\mid H_0)$ = $\alpha_1$,
\item otherwise reject $H_0$ for $H_2$.
\end{enumerate}
\end{definition}

A completely different approach {to testing inequality constraints, which provides full protection against
stating that the assumed inequality constraints hold when they are
violated, based on $\m{min}(\b z)$}, was proposed by \citet{Sasabuchi} and has been considered
by \citet{KatAgr}. This procedure has, however, very serious
drawbacks \citep[see][4.6]{DarFor98}.

Multiple comparison procedures are computationally simpler for two reasons: (i) they can be applied without the need to fit the inequality constrained model and (ii) they do not require computation of the probability weights of the chi bar squared distribution and exploit modern advances in the computation of probability integrals for the multivariate normal distribution. Critical values for the previous procedures rely on the computation of the probability $\Phi(\b a, \b b, \b \Sigma_0)$ that a multivariate normal $N(\b 0, \b \Sigma_0)$ lies in the rectangle $[ \b a, \b b ]$. For a survey of this problem in the context of multiple comparison procedures see \citet{GeBr2} chapter 3 and \citet{GeBr} chapter 6. Because  $\alpha_1$, $\alpha_2$ and $\alpha_{12}$ are probabilities defined under the null hypothesis, in order to compute  the critical values, the unknown  $\b \Sigma_0$ must be replaced by its estimate $\hat{\b \Sigma}_0$ under $H_0$. Alternatively the elements of $\b z$ can be divided by the  standard error estimated under $H_0$ and the critical values computed  by evaluating the previous multi-normal integrals  using the correlation matrix corresponding
to $\hat{\b \Sigma}_0$.
\section{Examples}
Table \ref{KAg} below was used by \citet{KatAgr} as an
instance of a context where the statistic $L_{01}$ can give
misleading evidence in favour of $H_1$ when we test $H_0$
(independence) against $H_1$ (all local log-odds ratios are non
negative).
\begin{table}[!h]
\centering
\caption{Treatment and Extent of trauma due to subarachnoid
hemorrhage}
\begin{tabular}{rrrrrr}
  \hline
 & Death & Veget & Major & Minor & Recov \\
  \hline
   Placebo & 59 & 25 &  46 &  48 &  32 \\
  Treated & 135 & 39 & 147 & 169 & 102 \\
   \hline
\end{tabular}
\label{KAg}
\end{table}
Here we have $L_{01}=7.89$ and  $L_{12}=1.75$, so, when we compare these statistics with the critical values of several tunable
LR procedures computed with exact weights and displayed in table \ref{cval} for $\alpha_1=0.02$ and
$\alpha_2=0.03$, we see that all procedures reject
$H_0$ in favour of $H_2$, except when the tuning parameter is
$\alpha_{12}=0$. Thus, in this case, a procedure with $\alpha_{12}=0.015$ seems to be a reasonable choice.  For comparison, we also apply a set of multiple comparisons
procedures.
\begin{table}[!th]
\centering {\small
\caption{Critical values for the data in Table \ref{KAg}}
\vskip1mm
\begin{tabular}{rrrrrrr}
  \hline \hline
$\alpha_{12}$ & 0 & 0.015 & 0.02 & 0.025 & 0.028 & 0.03 \\
\hline
 & \multicolumn{6}{c}{LR tunable procedures} \\
  \hline
$c_2$& 8.87 & 10.44 & 11.34 & 12.88 & 14.89 & Inf \\
  $c_1$ & 6.83 & 5.70 & 5.43 & 5.19 & 5.07 & 4.98 \\
  $ c_{12}$ & 8.87 & 1.46 & 1.16 & 0.95 & 0.86 & 0.81 \\
   \hline
 & \multicolumn{6}{c}{MC tunable procedures} \\
  \hline
$c_2$& 2.43 & 2.67 & 2.81 & 3.03 & 3.29 & Inf \\
  $c_1$ & 2.48 & 2.32 & 2.29 & 2.26 & 2.25 & 2.24 \\
  $c_{12}$ & 2.43 & 1.84 & 1.74 & 1.67 & 1.64 & 1.62 \\
   \hline
\end{tabular}
} \label{cval}
\end{table}
The minimum and maximum unconstrained estimates of the
log-odd ratios (studentized using the standard errors estimated
under $H_0$) are equal to -1.168 and 2.186 respectively. These
statistics are assessed against the critical values of several
tunable MC procedures displayed in Table \ref{cval}. All the procedures accept $H_0$, {a result which shows that the MC procedures tend to be more conservative}.

As a second example, we analyze the data in Table 3 from \citet{BaDaFo}; this is a $2\times 5\times 5$ contingency table where a sample of 2904 males were classified according to two age classes and 5 ordered categories of their own occupational prestige (OP) and that of their fathers. The first issue is whether association between father's and son's OP, {measured by log-odds ratio of a suitable type, is stronger when sons are older; so, this is equivalent to assume that log-odds ratios  increase with son age.}
If we use the global log-odds ratios, we have  $L_{01}=29.02$ and $L_{12}=0.29$ and  the  conclusion is that $H_0$ is rejected in favour of $H_1$ by any LR tunable procedure with $\alpha_1=0.02$ and $\alpha_2=0.03$. For instance, using again exact weights, with $\alpha_{12}=0$ the critical values are $c_1=25.197$, $c_2=c_{12}=11.120$, while with $\alpha_{12}=0.029$ they are  $c_1=22.117$, $c_2=19.731$, $c_{12}=1.394$; in both cases  $c_{12}$ is much larger than the observed value of $L_{12}$. Instead, if we measure the strength of association by the local log-odds ratios  we have $L_{01}=14.70$ and $L_{12}=12.55$ and, though we reject again $H_0$ towards $H_1$ with $\alpha_{12}=0$ because $c_1=12.214$, $c_2=c_{12}=23.820$, we reject towards $H_2$ if we set $\alpha_{12}=0.015$  because $c_1=10.746$, $c_2=26.161$, $c_{12}=10.785$. This example again shows that the LR tunable procedures allow us to submit $H_1$ to a severe scrutiny, while the standard procedure would often be too liberal relative to accepting the assumed ordering.

For the same data, it might be of interest to consider also the assumption that, when sons are younger, their OP is stochastically smaller than that of their fathers while the situation reverses when they are older; {the idea behind is that, due to welfare improvement, sons will, on the whole, be better off than their fathers, but also that OP improves with age}. This assumption compares the marginal distributions of fathers and sons conditionally on the age group of the sons, using logits of type global. Here $L_{01}=130.00$ and $L_{12}=5.12$ and, with $\alpha_1=0.02$, $\alpha_2=0.03$, $\alpha_{12}=0$, the critical values are $c_2=c_{12}=9.606$ and $c_1=14.211$, so $H_0$ is rejected in favour of  $H_1$  though there are substantial violations in the observed data. Instead, if we set $\alpha_{12}=0.015$, with the critical values $c_1=12.719$, $c_2=11.292$, $c_{12}=1.542$  $H_0$ must be rejected in favour of $H_2$.
This example shows that, for large values of $L_{01}$ and small values of $L_{12}$, the naive LR procedure  can give false evidence in favour of $H_1$ when $H_2$ is true and that this drawback is eliminated by the introduction of the tuning probability $\alpha_{12}.$

\section{Simulation study}
To evaluate the performance of the various procedures, we used a
targeted set of simulations concerning the problem of testing
independence against the assumption that all the log-odds ratios of
type local-local in $3\times3$ and $3\times 4$ contingency tables
are non negative. The sample size $n$ and the number of replications
$N$  were fixed to 10,000. To keep the context
simple, in all the simulations we set $\alpha_1=0.02$ and
$\alpha_2=0.03$ combined with a range of values for $\alpha_{12}$.
Initial simulations were used to check that all procedures achieved,
with high accuracy, the correct size under $H_0$, the hypothesis  of
independence; then, we considered various versions of $H_1$ and
$H_2$ by selecting specific sets of local odds ratios and
constructed the corresponding bivariate distributions having uniform
marginals.

The results are displayed in Tables 3-5 below. In all the tables,
each column corresponds to a different testing procedures. Within
the likelihood ratio (LR) approach, each procedure is identified by
the value of the tuning parameter $\alpha_{12}$, while, within the
multiple comparison (MC) approach, results for the Bennet procedure
are displayed in addition; for the values of $\alpha_1,\:\alpha_2$
chosen in the simulations, the rejection regions of the  naive MC procedure are very close to those of the Bennet procedure, thus
results are omitted because they either coincide or the difference
is smaller than $10^{-3}$. Tables are divided into sub-tables, each
corresponding to different experiments, identified by the set of
log-odds ratios used to generate the data. Each sub-table has three
rows giving the relative frequencies of: (i) acceptance of $H_0$,
(ii) rejections in favour of $H_1$ and (iii) rejections in favour of
$H_2$, for every procedure.
\subsection{Power under $H_1$}
All LR procedures seem to have high power even for moderate
violations of $H_0$ in the direction of $H_1$, like in Tables \ref{LR33}(a),\ref{LR34}(a);
with larger violations power gets very close to 1 as in Tables \ref{LR33}(b), \ref{LR34}(b).
Though we explored only the case where all the log-odds ratios were
equal, it seems reasonable to expect that the performance is
determined by the smallest positive value. The power decreases
slightly when $\alpha_{12}$ increases, this is consistent with the
fact that lager values of $\alpha_{12}$ provide more protection
against rejecting towards $H_1$ when $H_2$ is true. The power of MC
procedures is considerably smaller with a relatively large error
rate in the direction of $H_0$.
\subsection{Power under $H_2$}
Here the situation is much more complex because violations of $H_1$
can arise in many different directions and it is unlikely that a
procedure can perform best under all possible $H_2$ alternatives.
In the simulations we explored a limited range of possibilities
which, however, seem to suggest some general conclusions.

The first result which emerges clearly is that the naive LR procedure {cannot be recommended} because there are relevant versions of $H_2$ under which this procedure accepts $H_0$ with probability close to 1. This happens, for instance, when all log-odds ratios are negative, which means that all inequalities are violated, like in Table \ref{LR33}(g,h). The naive procedure has also a rather poor performance when the negative log-odds ratios dominate in number or in absolute value, like in Tables \ref{LR33}(i), \ref{LR34}(c). In all such cases the naive procedure tends to be terribly conservative. It is true that under few violations of $H_1$, like in Tables \ref{LR33}(c,d) and \ref{LR34}(d,e) the naive procedure is the best; {however, in these cases, most of the other tuned procedures have high power and the improvement produced by the naive procedure} is rather modest and certainly cannot compensates the very bad performance described above.
The performance of the tuned LR procedures seem to follow this
general pattern: when violations of $H_1$ in the direction of $H_2$
are few or the negative log-odds ratios are sufficiently small in
absolute value relatively to the positive ones, the procedure with
$\alpha_{12}=0$ tends to reject in the direction of $H_1$ with a
relatively large error rate which, however, can be reduced
dramatically by increasing $\alpha_{12}$, see for instance Tables
\ref{LR33}(c,d,j) and \ref{LR34}(d,e,f,g). Instead, when there is some kind of balance
between negative and positive log-odds ratios, like in Tables \ref{LR33}(e,k)
and \ref{LR34}(i), the LR procedures have a relatively large error rate in
the direction of $H_0$ which increases with $\alpha_{12}$.

On the whole, MC procedures perform substantially worst than the
corresponding LR procedure, with a few exceptions, typically when
the negative and positive log-odds ratios are is some kind of
balance: compare, for instance, {the corresponding entries in}
Table \ref{MC33}(e,k) and Table \ref{LR33}(e,k). MC procedures also do better when the
negative log-odds ratios are small or few in number but only for
$\alpha_{12}=0$; however, already at $\alpha_{12}=0.015$ the LR
procedure perform much better, like in Tables \ref{LR33},\ref{MC33}(c,d). Usually,
performance of MC procedures improve with $\alpha_{12}$, though they
remain too much inferior relative to the {corresponding} LR
procedures. The MC procedures are again substantially inferior when
all log-odds ratios are negative, like in Tables \ref{MC33}(g,h). On the
whole, Bennet's procedure performs like an MC procedure with
$\alpha_{12}$ around 0.010; this, however, may depend on the
specific values of $\alpha_1,\:\alpha_2$ used in the simulation.
Note also that the naive MC procedure does not do as badly as in the
LR context.

\newpage
\begin{table}[!h]
\centering \caption{Likelihood ratio procedures in $3\times 3$
tables} {\small
\begin{tabular}{rrrrrrr}
  \hline \hline
 & \multicolumn{6}{c}{$\alpha_{12}$}\\
 & 0.000 & 0.015 & 0.020 & 0.025 & 0.028 & 0.030 \\
  \hline\hline
 & \multicolumn{6}{c}{(a), $H_1$: 0.08, 0.08, 0.08, 0.08}\\ \hline
H0 & 0.002 & 0.001 & 0.001 & 0.000 & 0.000 & 0.000 \\
H1 & 0.998 & 0.975 & 0.965 & 0.952 & 0.945 & 0.941 \\
H2 & 0.000 & 0.025 & 0.035 & 0.048 & 0.055 & 0.059 \\
 \hline
 & \multicolumn{6}{c}{(b), $H_1$: 0.15, 0.15, 0.15, 0.15}\\ \hline
H0 & 0.000 & 0.000 & 0.000 & 0.000 & 0.000 & 0.000 \\
H1 & 1.000 & 1.000 & 0.999 & 0.999 & 0.999 & 0.999 \\
H2 & 0.000 & 0.000 & 0.001 & 0.001 & 0.001 & 0.001 \\
   \hline
 & \multicolumn{6}{c}{(c), $H_2$: 0.08, 0.08, -0.08, 0.08}\\ \hline
H0 & 0.172 & 0.121 & 0.111 & 0.104 & 0.100 & 0.096 \\
H1 & 0.769 & 0.351 & 0.305 & 0.271 & 0.253 & 0.245 \\
H2 & 0.059 & 0.528 & 0.584 & 0.626 & 0.647 & 0.658 \\
   \hline
 & \multicolumn{6}{c}{(d), $H_2$: 0.15, 0.15, -0.15, 0.15}\\ \hline
H0 & 0.000 & 0.000 & 0.000 & 0.000 & 0.000 & 0.000 \\
H1 & 0.690 & 0.092 & 0.072 & 0.060 & 0.054 & 0.051 \\
H2 & 0.310 & 0.908 & 0.928 & 0.940 & 0.946 & 0.949 \\
   \hline
 & \multicolumn{6}{c}{(e), $H_2$: -0.08, 0.08, 0.08, -0.08}\\ \hline
H0 & 0.785 & 0.828 & 0.845 & 0.867 & 0.889 & 0.919 \\
H1 & 0.034 & 0.010 & 0.009 & 0.009 & 0.008 & 0.008 \\
H2 & 0.181 & 0.162 & 0.146 & 0.124 & 0.103 & 0.074 \\
   \hline
 & \multicolumn{6}{c}{(f), $H_2$: -0.15, 0.15, 0.15, -0.15}\\ \hline
H0 & 0.305 & 0.380 & 0.424 & 0.495 & 0.581 & 0.833 \\
  H1 & 0.031 & 0.002 & 0.001 & 0.001 & 0.001 & 0.001 \\
  H2 & 0.664 & 0.619 & 0.575 & 0.503 & 0.419 & 0.167 \\
   \hline
 & \multicolumn{6}{c}{(g), $H_2$: -0.08, -0.08, -0.08, -0.08}\\ \hline
H0 & 0.003 & 0.007 & 0.009 & 0.016 & 0.033 & 1.000 \\
H1 & 0.000 & 0.000 & 0.000 & 0.000 & 0.000 & 0.000 \\
H2 & 0.997 & 0.993 & 0.991 & 0.984 & 0.967 & 0.000 \\
   \hline
 & \multicolumn{6}{c}{(h), $H_2$: -0.15, -0.15, -0.15, -0.15}\\ \hline
H0 & 0.000 & 0.000 & 0.000 & 0.000 & 0.000 & 1.000 \\
H1 & 0.000 & 0.000 & 0.000 & 0.000 & 0.000 & 0.000 \\
H2 & 1.000 & 1.000 & 1.000 & 1.000 & 1.000 & 0.000 \\
   \hline
 & \multicolumn{6}{c}{(i), $H_2$: -0.15, 0.04, 0.04, -0.15}\\ \hline
H0 & 0.047 & 0.080 & 0.105 & 0.154 & 0.230 & 1.000 \\
H1 & 0.000 & 0.000 & 0.000 & 0.000 & 0.000 & 0.000 \\
H2 & 0.953 & 0.920 & 0.895 & 0.846 & 0.770 & 0.000 \\
   \hline
 & \multicolumn{6}{c}{(j), $H_2$: -0.04, 0.15, 0.15, -0.04}\\ \hline
H0 & 0.044 & 0.028 & 0.024 & 0.022 & 0.020 & 0.020 \\
H1 & 0.925 & 0.468 & 0.412 & 0.373 & 0.353 & 0.340 \\
H2 & 0.031 & 0.504 & 0.563 & 0.605 & 0.627 & 0.641 \\
   \hline
 & \multicolumn{6}{c}{(k), $H_2$: -0.12, 0.08, 0.08, -0.12}\\ \hline
0 & 0.469 & 0.574 & 0.626 & 0.708 & 0.795 & 0.993 \\
H1 & 0.001 & 0.000 & 0.000 & 0.000 & 0.000 & 0.000 \\
H2 & 0.530 & 0.426 & 0.374 & 0.292 & 0.205 & 0.007 \\

   \hline   \hline
\end{tabular}
}\label{LR33}
\end{table}

\newpage
\begin{table}[!h]
\centering \caption{Likelihood ratio procedures in $3\times 4$
tables} {\small
\begin{tabular}{rrrrrrr}
  \hline \hline
 & \multicolumn{6}{c}{$\alpha_{12}$}\\
 & 0.000 & 0.015 & 0.020 & 0.025 & 0.028 & 0.030 \\
  \hline\hline
  & \multicolumn{6}{c}{(a), $H_2$: 0.06,0.06,0.06,0.06,0.06,0.06}\\ \hline
  H0 & 0.001 & 0.000 & 0.000 & 0.000 & 0.000 & 0.000 \\
  H1 & 0.999 & 0.960 & 0.946 & 0.930 & 0.921 & 0.916 \\
  H2 & 0.000 & 0.040 & 0.054 & 0.069 & 0.078 & 0.084 \\
  \hline
  & \multicolumn{6}{c}{(b), $H_2$: 0.12,0.12,0.12,0.12,0.12,0.12}\\ \hline
 0 & 0.000 & 0.000 & 0.000 & 0.000 & 0.000 & 0.000 \\
  H1 & 1.000 & 0.998 & 0.997 & 0.995 & 0.994 & 0.993 \\
  H2 & 0.000 & 0.002 & 0.003 & 0.005 & 0.006 & 0.007 \\
  \hline
  & \multicolumn{6}{c}{(c), $H_2$: -0.17,0.15,0.15,-0.17,-0.17,0.15}\\ \hline
  H0 & 0.172 & 0.248 & 0.292 & 0.371 & 0.474 & 0.942 \\
  H1 & 0.005 & 0.000 & 0.000 & 0.000 & 0.000 & 0.000 \\
  H2 & 0.823 & 0.751 & 0.708 & 0.629 & 0.526 & 0.058 \\
 \hline
 & \multicolumn{6}{c}{(d), $H_2$: 0.12,0.12,-0.12,0.12,0.12,0.12}\\ \hline
H0 & 0.000 & 0.000 & 0.000 & 0.000 & 0.000 & 0.000 \\
  H1 & 0.952 & 0.483 & 0.428 & 0.390 & 0.367 & 0.356 \\

  H2 & 0.048 & 0.517 & 0.572 & 0.610 & 0.633 & 0.644 \\
   \hline
 & \multicolumn{6}{c}{(e), $H_2$: 0.15,0.15,-0.15,0.15,0.15,0.15}\\ \hline
H0 & 0.000 & 0.000 & 0.000 & 0.000 & 0.000 & 0.000 \\
  H1 & 0.902 & 0.331 & 0.285 & 0.251 & 0.235 & 0.227 \\
  H2 & 0.098 & 0.669 & 0.715 & 0.749 & 0.765 & 0.773 \\
   \hline
& \multicolumn{6}{c}{(f), $H_2$: -0.6,0.15,0.15,-0.6,-0.6,0.15}\\
   \hline
H0 & 0.021 & 0.012 & 0.010 & 0.009 & 0.008 & 0.008 \\
H1 & 0.932 & 0.433 & 0.377 & 0.337 & 0.316 & 0.306 \\
H2 & 0.047 & 0.555 & 0.613 & 0.654 & 0.676 & 0.685 \\
 \hline
 & \multicolumn{6}{c}{(g), $H_2$: 0.15,0.15,-0.15,0.15,0.15,0.15}\\
  \hline
H0 & 0.000 & 0.000 & 0.000 & 0.000 & 0.000 & 0.000 \\
  H1 & 0.902 & 0.331 & 0.285 & 0.251 & 0.235 & 0.227 \\
  H2 & 0.098 & 0.669 & 0.715 & 0.749 & 0.765 & 0.773 \\
   \hline
& \multicolumn{6}{c}{(i), $H_2$: -0.15,0.15,0.15,-0.15,-0.15,0.15}\\
 \hline
H0 & 0.316 & 0.386 & 0.432 & 0.503 & 0.579 & 0.814 \\
  H1 & 0.037 & 0.003 & 0.002 & 0.002 & 0.001 & 0.001 \\
  H2 & 0.647 & 0.611 & 0.566 & 0.496 & 0.420 & 0.184 \\
   \hline   \hline
\end{tabular}
}\label{LR34}
\end{table}

\newpage
\begin{table}[!h]
\centering \caption{Multiple comparisons procedures in $3\times 3$ tables
} {\small
\begin{tabular}{rrrrrrrr}
  \hline \hline
 & & \multicolumn{6}{c}{$\alpha_{12}$}\\
 & Bennet & 0.000 & 0.015 & 0.020 & 0.025 & 0.028 & 0.030 \\
  \hline\hline
 &\multicolumn{7}{c}{(a),$H_1$: 0.08, 0.08, 0.08, 0.08}\\ \hline
H0 & 0.554 & 0.573 & 0.464 & 0.438 & 0.417 & 0.406 & 0.401  \\
  H1 & 0.445 & 0.427 & 0.533 & 0.558 & 0.578 & 0.589 & 0.593 \\
  H2 & 0.001 & 0.001 & 0.004 & 0.004 & 0.005 & 0.005 & 0.006 \\
   \hline
 &\multicolumn{7}{c}{(b), $H_1$: 0.15, 0.15, 0.15, 0.15}\\ \hline
H0 & 0.008 & 0.010 & 0.003 & 0.002 & 0.001 & 0.001 & 0.001  \\
  H1 & 0.992 & 0.990 & 0.997 & 0.999 & 0.999 & 0.999 & 0.999  \\
  H2 & 0.000 & 0.000 & 0.000 & 0.000 & 0.000 & 0.000 & 0.000 \\
   \hline
 &\multicolumn{7}{c}{(c), $H_2$: 0.08, 0.08, -0.08, 0.08}\\ \hline
H0 & 0.598 & 0.608 & 0.547 & 0.533 & 0.520 & 0.516 & 0.514  \\
  H1 & 0.249 & 0.255 & 0.234 & 0.230 & 0.228 & 0.228 & 0.229  \\
  H2 & 0.154 & 0.137 & 0.219 & 0.238 & 0.251 & 0.256 & 0.257 \\
   \hline
 &\multicolumn{7}{c}{(d), $H_2$: 0.15, 0.15, -0.15, 0.15}\\ \hline
H0 & 0.044 & 0.047 & 0.025 & 0.021 & 0.019 & 0.018 & 0.017 \\
  H1 & 0.373 & 0.424 & 0.229 & 0.201 & 0.185 & 0.177 & 0.172  \\
  H2 & 0.583 & 0.529 & 0.746 & 0.778 & 0.797 & 0.805 & 0.810  \\
   \hline
 &\multicolumn{7}{c}{(e), $H_2$: -0.08, 0.08, 0.08, -0.08}\\ \hline
H0 & 0.659 & 0.658 & 0.664 & 0.670 & 0.676 & 0.682 & 0.685  \\
  H1 & 0.082 & 0.097 & 0.048 & 0.041 & 0.037 & 0.035 & 0.034  \\
  H2 & 0.259 & 0.245 & 0.287 & 0.289 & 0.287 & 0.283 & 0.281 \\
   \hline
 &\multicolumn{7}{c}{(f), $H_2$: -0.15, 0.15, 0.15, -0.15}\\ \hline
H0 & 0.165 & 0.164 & 0.168 & 0.174 & 0.182 & 0.189 & 0.195  \\
  H1 & 0.071 & 0.100 & 0.022 & 0.016 & 0.012 & 0.011 & 0.010 \\
  H2 & 0.764 & 0.737 & 0.810 & 0.810 & 0.806 & 0.800 & 0.794  \\
   \hline
 &\multicolumn{7}{c}{(g), $H_2$: -0.08, -0.08, -0.08, -0.08}\\ \hline
H0 & 0.541 & 0.524 & 0.675 & 0.742 & 0.827 & 0.902 & 0.999  \\
  H1 & 0.000 & 0.000 & 0.000 & 0.000 & 0.000 & 0.000 & 0.000  \\
  H2 & 0.459 & 0.476 & 0.325 & 0.258 & 0.173 & 0.098 & 0.001  \\
   \hline
 &\multicolumn{7}{c}{(h), $H_2$: -0.15, -0.15, -0.15, -0.15}\\ \hline
H0 & 0.009 & 0.008 & 0.038 & 0.067 & 0.153 & 0.313 & 1.000  \\
  H1 & 0.000 & 0.000 & 0.000 & 0.000 & 0.000 & 0.000 & 0.000  \\
  H2 & 0.991 & 0.992 & 0.962 & 0.933 & 0.847 & 0.687 & 0.000  \\
   \hline
 &\multicolumn{7}{c}{(i), $H_2$: -0.15, 0.04, 0.04, -0.15}\\ \hline
H0 & 0.274 & 0.262 & 0.360 & 0.416 & 0.515 & 0.621 & 0.891  \\
  H1 & 0.000 & 0.000 & 0.000 & 0.000 & 0.000 & 0.000 & 0.000 \\
  H2 & 0.726 & 0.737 & 0.640 & 0.584 & 0.485 & 0.379 & 0.109  \\
 \hline
&\multicolumn{7}{c}{(j), $H_2$: -0.04, 0.15, 0.15, -0.04}\\ \hline
H0 & 0.265 & 0.277 & 0.219 & 0.206 & 0.195 & 0.189 & 0.186  \\
  H1 & 0.632 & 0.649 & 0.559 & 0.536 & 0.520 & 0.514 & 0.510  \\
  H2 & 0.104 & 0.074 & 0.222 & 0.258 & 0.285 & 0.297 & 0.304  \\
   \hline
 &\multicolumn{7}{c}{(k), $H_2$: -0.12, 0.08, 0.08, -0.12}\\ \hline
H0 & 0.460 & 0.453 & 0.520 & 0.554 & 0.600 & 0.641 & 0.695  \\
  H1 & 0.020 & 0.028 & 0.006 & 0.005 & 0.005 & 0.004 & 0.004  \\
  H2 & 0.520 & 0.520 & 0.473 & 0.441 & 0.396 & 0.355 & 0.302 \\
   \hline \hline
\end{tabular}
}\label{MC33}
\end{table}

\newpage
\section*{Appendix}
{\bf Computation of probability weights}

In order to compute  the weights $w_i(\cg C,\b V)$, it may be
useful to summarize the geometry of the projection of a
random vector $\b y \sim \cg N(\b 0,\b V)$ onto a convex cone $ \cg
C=\{\b \eta:\b D\b\eta\geq \b 0\}$, where  $\b D$ is a $k\times
(t-1)$ matrix of rank $k$. If $\b H$ is the left component of the
Cholesky decomposition of the positive definite matrix  $\b\Omega$ =
$\b D \b V \b D\tr$ then $\b z$ = $\b H^{-1} \b D\b y$ $\sim$ $N(\b
0,\b I_k)$, the transformation $\b \lambda$ = $\b H^{-1} \b D\b
\eta$ defines the cone $ \cg C^* =\{\b \lambda:\b H \b\lambda\geq \b
0 \}$,  $\cg C \in \Re^k$, and it holds that: $\min_{\bl
D\bl\eta\geq \bl 0}(\b y-\b\eta)\tr\b V^{-1}(\b y-\b\eta)=\min_{\bl
H\bl\lambda\geq \bl 0}(\b z-\b\lambda)\tr(\b z-\b\lambda)$. The
cone  $\cg C^*$ may also be defined by its generating vectors which
are the columns of  $\b U = \b H^{-1}$: a vector $\b z$ belongs to
$\cg C^*$ if $\b z=\b U\b u$ where $\b u\geq \b 0$. In a similar way
the dual cone $\cg C^{*0}$ is generated by the columns of $\b W=-\b
H\tr$ and note that $\b U\tr\b W$ = $-\b I$.

Within the euclidian metric, $\Re^k$ can be partitioned into $2^k$
convex cones as follows: let $\cg J$ be the collection of all
possible subsets of $(1, \dots, k)$, including the empty set and the
whole set. For any pair $\b i,\:\b j\in\cg J$, $\b i\cup\b j$ = $(1,
\dots, k)$, let $\left(\b U_{\bl i},\: \b W_{\bl j}\right)$, be the
matrix whose columns are, respectively, the columns of $\b U$ with
index in $\b i$ and the columns of $\b W$ with index in $\b j$; the
columns of this matrix generate the convex cone $\cg C^*(\b i)$
whose elements, when projected onto $\cg C^*$, belong to the face
generated by the columns of $\b U_{\bl i}$, this face is itself a
convex cone of dimension equal to the cardinality $|\b i|$ of $\b
i$. Thus
$$
w_{i+q}(\cg C,\b V) =w_i(\cg C^*,\b I) =\sum_{\mid \bl i\mid =i}
P[\b z\in \cg C^*(\b i)], \quad i=0,1,...,k
$$
where $q$ is the dimension of $\cg L_0$.

To compute $P[\b z\in \cg C^*(\b i)]$ note that $\b z\in \cg C^*(\b
i)$ if and only if $\b t$ = $\left(\b U_{\bl i},\: \b W_{\bl
j}\right)^{-1}\b z\geq \b 0$, in other words, the linear
transformation above reduces $\cg C^*(\b i)$ into the positive
orthant for the multivariate normal random variable $\b t$; thus,
{to compute $P[\b t\in \cg R^{k+}]$, the only quantity we need is
$Var(\b t)$ = $\b\Omega$. Let $\b\Psi$ = $\b W\tr\b W$ and $\b\Phi$
= $(\b U\tr\b U)$ and note that $\b\Psi$ = $\b D\b V\b D\tr$ =
$\b\Phi^{-1}$. It can be shown that $\b\Omega$ is block diagonal
with elements given by $(\b\Phi_{\bl i\bl i})^{-1}$ and
$(\b\Psi_{\bl j\bl j})^{-1}$,} which are related by the well known
formulas for the inverse of a partitioned matrix:
\begin{eqnarray*}
(\b\Phi_{\bl i\bl i})^{-1} &=& \b\Psi_{\bl i\bl i}-\b\Psi_{\bl i\bl
j}(\b\Psi_{\bl j\bl j})^{-1}\b\Psi_{\bl j\bl i}\\
(\b\Psi_{\bl j\bl j})^{-1} &=& \b\Phi_{\bl j\bl j}-\b\Phi_{\bl j\bl
i}(\b\Phi_{\bl i\bl i})^{-1}\b\Phi_{\bl i\bl j}.
\end{eqnarray*}
So, if $\mid\b i\mid \leq \mid j\mid$, it is convenient to
compute $(\b\Phi_{\bl i\bl i})^{-1}$ directly and $(\b\Psi_{\bl j\bl
j})^{-1}$ from the second expression above, instead, when $\mid\b
i\mid > \mid j\mid$, compute $(\b\Psi_{\bl j\bl j})^{-1}$  directly
and  $(\b\Phi_{\bl i\bl i})^{-1}$ from the first expression above.
{In any case, because $\b\Omega$ is block diagonal,
$P[\b t\in \cg R^{k+}]$ factorizes into the product of two lower dimensional integrals.}

Because Proposition 3.6.1(3) in \citet{SiSe} says that the weights
with index $j$ even or odd sum to 0.5, we may avoid computing the
two weight which correspond to the largest number of side cones;
these correspond to $(k/2-1,\:k/2)$ when $k$ is even and to
$((k-1)/2,\:(k+1)/2)$ when $k$ is odd.

\section*{References}
\bibliographystyle{plainnat}
\bibliography{ChiBar}

\end{document}